# Power System Simulation Using the Differential Transformation Method

Yang Liu, *Student Member, IEEE*, and Kai Sun, *Senior Member, IEEE*

*Abstract*— **This paper proposes a new semi-analytical approach for online time-domain power system simulation. The approach applies the differential transformation method (DTM) to the power system differential equation model to offline derive a semi-analytical solution (SAS) having symbolic variables about time, the initial state and system conditions. When simulation is online needed for a contingency under the current system condition, the SAS can be evaluated in real time to generate simulation results. Compared to the Adomian decomposition method in obtaining a power system SAS, an SAS derived by the DTM adopts a recursive form to avoid generating and storing its complete symbolic expression, which makes both derivation and evaluation of the SAS more efficient especially for multi-machine power systems. The optimal order of a DTM-based SAS is studied for the best time performance of simulation. The paper also designs a parallel computing strategy for power system simulation using the DTM-based SAS. Tests on the IEEE 10-machine 39-bus system demonstrate significant speedup of simulation using the proposed approach compared with the Runge-Kutta method.**

## I. INTRODUCTION

**P**OWER system transient stability simulation is of critical importance for utilities to assess the dynamic security by solving the initial value problem (IVP) of nonlinear differential equations (DEs) with a given contingency occurring under a specific operating condition. Numerical integration methods, including explicit and implicit methods, are commonly used in commercial software packages with a small enough integration step of typically one to few milliseconds to meet accuracy requirements. In recent years, the power systems have been pushed to be operated closer to their stability limits due to the fast growth in electricity demands but a relatively slow construction of new transmission infrastructure. To identify any insecure contingency before it happens, time domain simulation is expected to be transitioned from offline or day-ahead studies to the real-time operation environment. The power industry and the research community are seeking next-generation simulation tools which are more powerful for power system dynamic security assessment in a faster than real-time manner.

One way to speed up transient stability simulation using traditional numerical integration methods is parallel computing. With the fast development of high performance computers (HPCs), a variety of parallel power system simulation methods have been proposed by means of decomposing the system model or computation tasks in simulation. Paper [1] decomposes a system model into three linear subsystems and ensures simulation accuracy by adaptive updates on linear subsystems. Paper [2] proposes a two-stage parallel waveform relaxation method for parallel simulation, which adopts epsilon decomposition to partition a large-scale power system model into several subsystems. A Schur-complement based network decomposition method is proposed in [3]. The parareal in time method in [4] and [5] adopts temporal decomposition of the simulation period into many intervals, conducts parallel simulations on individual intervals using a fine solver, and connects their results by a high-level coarse solver after a few iterations. Recently, paper [6] develops a practical framework to parallelize computation tasks of a single dynamic simulation in commercial software and paper [7] designs a massively parallel computational platform for efficient dynamic security assessment. All the above simulation approaches are based on numerical integration methods, which originally were not developed for a parallel computing environment, and hence can only be parallelized to a limit extent. Thus, emerging, alternative simulation methods out of the numerical integration framework are attracting more interests in research.

As a new paradigm for fast transient stability simulation, a semi-analytical approach is recently proposed and investigated in [8]-[13]. Its basic idea is to partition the computation towards solution of the IVP of power system DEs into two stages. The first stage off-line derives a semi-analytical solution (SAS) that expresses each state variable as an explicit function of symbolic variables including time, the initial state and parameters on system conditions. The SAS is an approximate but analytical solution being accurate for a certain time window whose length depends on the inherent nonlinearity and order of the SAS expression. The second stage on-line evaluates the SAS over consecutive time windows to make up the desired simulation period by substituting values for those symbolic variables according to a given contingency and the real-time system condition. Three methods, i.e. Taylor Expansion (TE) [8], Adomian Decomposition method (ADM) [9]-[11], and Pade Approximants (PA) [12], have been applied to power system models to derive SAS's respectively in the forms of

This work was supported by grant ECCS-1610025.

Y. Liu and K. Sun are with the Department of Electrical Engineering and Computer Science, University of Tennessee, Knoxville, TN 37996 USA (e-mail: yliu161@vols.utk.edu, kaisun@utk.edu).

polynomials of time, sinusoidal or polynomial functions of time, and fractional functions of time. Ref. [13] further integrates ADM-based SAS's with numerical integration methods for a hybrid strategy and tests it on realistic large system models.

The speed of simulation using an SAS relies on the time window to maintain its accuracy, which further relies on its order, i.e. the number of summated terms, of its expression. However, when applied to a multi-machine power system, an SAS derived from the above three methods has considerably increased complexity especially after its order exceeds 3 because its expression includes all symbolic variables defined.

In this paper, a new approach for deriving a power system SAS in the form of polynomials is proposed based on the differential transformation method (DTM) [16]. A DTM-based SAS not only has the advantages of an SAS thanks to its analytical nature but also inherits some good properties from the DTM, in which the symbolic coefficients of different polynomial terms can be calculated recursively from low orders to high orders and stored as a set of concise and general recursive formulas that are applicable to any high order. By contrast, an ADM-based SAS has its coefficients derived individually by symbolic computations, which make the complete SAS expression very complex. By using a recursive form, a DTM-based SAS allows much more terms to be evaluated than an ADM-based SAS, so it can keep accuracy for much longer time to speed up simulation.

The paper has the following contributions: (1) it derives a DTM-based SAS for a multi-machine power system, analyzes its convergence and proposes an SAS-based transient stability simulation scheme by evaluating the SAS over multi-time windows; (2) the paper studies the optimal order of an SAS to be evaluated for fastest simulation, considering that an higher-order SAS has a longer time window of accuracy and hence needs fewer evaluations but it has a more complex expression to be evaluated; (3) the paper also studies the use of parallel computers for online evaluation of a DTM-based SAS to discover the maximum speeding-up of transient stability simulation by using the SAS.

The rest of the paper is organized as follows. Section II first briefly introduces the DTM and its formulation and then derives the DTM-based SAS of a multi-machine power system with detailed generator models. Section III presents the online scheme for transient stability simulation using a DTM-based SAS and the implementation strategy on parallel computers. Section IV uses the IEEE 39-bus system to compare the accuracy and time performance of the proposed semi-analytical approach with the Runge-Kutta-4th-order (RK4) method. Finally, conclusions are drawn in section V.

## II. THE DTM BASED SAS ON MULTI-MACHINE SYSTEMS

### A. Introduction of the DTM

The DTM theory is originally established in [16] to derive approximate solutions of nonlinear DEs by means of the differential transformation (DT) defined below. It is then developed by researchers in the fields of applied mathematics and physics to obtain SAS's of various nonlinear dynamic systems such as the Van der Pol oscillator, Duffing equations and fractional order systems [17]-[24]. In existing literature, the DTM is mainly applied to small systems described by low-order DEs and its capability has not been examined for real-life complex network systems like power systems modeled by high-order nonlinear DEs.

**Definition 1:** Consider the function $x = x(t)$ of a real continuous variable $t$. The differential transformation (DT) of $x(t)$ is defined in (1), and the inverse differential transformation (IDT) of $X(k)$ is defined in (2), where $k$ is the DT order.

$$X(k) = \frac{1}{k!}\left[\frac{d^k x(t)}{dt^k}\right]_{t=0} \tag{1}$$

$$x(t) = \sum_{k=0}^{\infty} X(k) \, t^k \tag{2}$$

The notation $x(t)$ is the original function in the time domain and $X(k) = \mathrm{DT}(x(t))$ is in the domain about order $k$. In this paper, lowercase and capital letters are respectively used for original functions and their DTs for differentiation.

The following Propositions are necessary to derive the SAS in this paper. Proposition 1 gives six basic propositions of the DTM. Corollary 1 is about the DT of a constant, which is easily deduced with $n=0$ in Proposition 1-(e).

**Proposition 1:** Denote $x(t)$, $y(t)$ and $z(t)$ as the original functions and $X(k)$, $Y(k)$ and $Z(k)$ as their DTs, respectively. The following propositions (a) - (f) hold for all continuous functions, where $c$ is a constant, $n$ is a nonnegative integer and $\boldsymbol{\delta}$ is the Dirac delta function

(a) $X(0) = x(0)$.

(b) If $y(t) = cx(t)$, then $Y(k) = cX(k)$.

(c) If $z(t) = x(t) \pm y(t)$, then $Z(k) = X(k) \pm Y(k)$.

(d) If $z(t) = x(t)y(t)$, then $Z(k) = \sum_{m=0}^{k} X(m)Y(k-m)$.

(e) If $z(t) = x^n(t)$, then $Z(k) = \boldsymbol{\delta}(k-n) = \begin{cases} 1, & k = n \\ 0, & k \neq n \end{cases}$.

(f) If $y(t) = \frac{dx(t)}{dt}$, then $Y(k) = (k+1)Z(k+1)$.

**Corollary 1:** Suppose $z(t) = c$ is a constant, then its DT is

$$Z(k) = c\boldsymbol{\delta}(k) = \begin{cases} c, & k = 0 \\ 0, & k \neq 0 \end{cases} \tag{3}$$

### B. Power System Model

Consider an $N$-machine power system modeled by (4)-(5), where for generator $i$, state variables $\delta_i$, $\omega_i$, $e'_{qi}$ and $e'_{di}$ are respectively the rotor angle and speed and $q$-axis and $d$-axis transient voltages, $P_{ei}$, $i_{di}$ and $i_{qi}$ are the electrical power and $d$-axis and $q$-axis stator currents, $e_{di}$ and $e_{qi}$ are the $d$-axis and $q$-axis terminal voltages, $i_{xi}$ and $i_{yi}$ are the $x$-axis and $y$-axis terminal currents, $e_{xi}$ and $e_{yi}$ are the $x$-axis and $y$-axis terminal voltages respectively, $H_i$ is the inertia, $T'_{d0i}$ and $T'_{q0i}$ are the open circuit transient time constants in $d$-axis and $q$-axis, $P_{mi}$ is the mechanical power, $e_{fdi}$ is field voltage, $x_{di}$, $x_{qi}$, $x'_{di}$ and $x'_{qi}$ are the $d$-axis and $q$-axis synchronous and transient reactances, $R_{ai}$ is the armature resistance and $D_i$ is the damping constant. Finally, $\omega_s = 2\pi \times 60$ is the nominal frequency and $\mathbf{Y}$ is the reduced network admittance matrix.

$$\begin{cases} \dot{\delta}_i = \omega_s \omega_i \\ \dot{\omega}_i = \frac{1}{2H_i}(P_{mi} - P_{ei} - D_i\omega_i) \\ \dot{e}'_{qi} = \frac{1}{T'_{d0i}}\left[e_{fdi} - e'_{qi} - (x_{di} - x'_{di})i_{di}\right] \\ \dot{e}'_{di} = \frac{1}{T'_{q0i}}\left[-e'_{di} - (x_{qi} - x'_{qi})i_{qi}\right] \end{cases} \quad i=1, ..., N \quad (4)$$

$P_{ei} = e_{di}I_{di} + e_{qi}I_{qi}$ where

$$\begin{bmatrix} i_{di} \\ i_{qi} \end{bmatrix} = \begin{bmatrix} \sin\delta_i & \cos\delta_i \\ -\cos\delta_i & \sin\delta_i \end{bmatrix}\begin{bmatrix} i_{xi} \\ i_{yi} \end{bmatrix}, \begin{bmatrix} e'_{di} \\ e_{qi} \end{bmatrix} = \begin{bmatrix} e'_{di} \\ e'_{qi} \end{bmatrix} - \begin{bmatrix} R_{ai} & -x'_{qi} \\ x'_{di} & R_{ai} \end{bmatrix}\begin{bmatrix} i_{di} \\ i_{qi} \end{bmatrix},$$

$$\begin{bmatrix} i_{x1} + ji_{y1} \\ \vdots \\ i_{xN} + ji_{yN} \end{bmatrix} = \mathbf{Y}\begin{bmatrix} e_{x1} + je_{y1} \\ \vdots \\ e_{xN} + je_{yN} \end{bmatrix}, \begin{bmatrix} e_{xi} \\ e_{yi} \end{bmatrix} = \begin{bmatrix} \sin\delta_i & \cos\delta_i \\ -\cos\delta_i & \sin\delta_i \end{bmatrix}\begin{bmatrix} e'_{qi} \\ e'_{di} \end{bmatrix}$$ (5)

### C. DTM-based SAS for a Multi-Machine System

Let $\boldsymbol{\varphi}(t)$ denote state vector $[\delta_i, \omega_i, e'_{qi}, e'_{di}]^\mathrm{T}$ and $\boldsymbol{\psi}(t)$ denote $[i_{di}, i_{qi}, P_{ei}, e_{xi}, e_{yi}, i_{xi}, i_{yi}, e_{di}, e_{qi}]^\mathrm{T}$, a vector of algebraic variables. Then the model (4)-(5) can be written as (6) where algebraic variables can be eliminated by substituting them into the DEs.

$$\dot{\boldsymbol{\varphi}} = \mathbf{f}(\boldsymbol{\varphi}, \boldsymbol{\psi}), \quad \boldsymbol{\psi} = \mathbf{g}(\boldsymbol{\varphi}) \quad (6)$$

In existing literature, the exact analytical solution does not exist even for an SMIB (single-machine-infinite-bus) model. The main difficulty lies in the nonlinear sine and cosine functions in $\boldsymbol{\psi} = \mathbf{g}(\boldsymbol{\varphi})$. Due to this obstacle, the DEs are solved by numerical integration in power system simulation. However, they can be solved analytically using the DTM as follows.

#### 1) Basic idea of a DTM-based SAS

Although it is impossible to obtain an exact analytical solution, trajectories of variables in $\boldsymbol{\varphi}(t)$ and $\boldsymbol{\psi}(t)$ can be approximated by polynomial functions of time in (7), where $K$ is their order and $\boldsymbol{\Phi}(k)$ and $\boldsymbol{\Psi}(k)$ are coefficients, which are functions of unknown symbolic variables about the initial state and system conditions.

$$\boldsymbol{\varphi}(t) = \sum_{k=0}^{K} \boldsymbol{\Phi}(k)\, t^k, \quad \boldsymbol{\psi}(t) = \sum_{k=0}^{K} \boldsymbol{\Psi}(k)\, t^k \quad (7)$$

Now the key to obtain trajectories $\boldsymbol{\varphi}(t)$ and $\boldsymbol{\psi}(t)$ is calculation of coefficients $\boldsymbol{\Phi}(k)$ and $\boldsymbol{\Psi}(k)$ ($k = 0, 1, \cdots K$). From (2), it is obvious that $\boldsymbol{\Phi}(k) = \mathrm{DT}(\boldsymbol{\varphi}(t))$ and $\boldsymbol{\Psi}(k) = \mathrm{DT}(\boldsymbol{\psi}(t))$. Thus, by applying DT to both sides of (7), we can solve $\boldsymbol{\Phi}(k)$ and $\boldsymbol{\Psi}(k)$ using Propositions 1. The DTM based SAS derivation can be summarized as the following steps.

**Step 1**: Apply the DT to functions on both sides of (6) to obtain (8), where $\boldsymbol{\Phi} = \mathrm{DT}(\boldsymbol{\varphi})$, $\boldsymbol{\Psi} = \mathrm{DT}(\boldsymbol{\psi})$, $\mathbf{F} = \mathrm{DT}(\mathbf{f})$ and $\mathbf{G} = \mathrm{DT}(\mathbf{g})$ based on Propositions.

$$\begin{aligned} \boldsymbol{\Phi}(k) &= \mathbf{F}\big(\boldsymbol{\Phi}(k-1), \boldsymbol{\Psi}(k-1)\big), k = 1 \cdots K \\ \boldsymbol{\Psi}(k) &= \mathbf{G}\big(\boldsymbol{\Phi}(0) .. \boldsymbol{\Phi}(k)\big), k = 1 \cdots K \end{aligned} \quad (8)$$

**Step 2**: Set a desired order $K$. Then staring from the initial state, all the coefficients from order 0 to order $k$ are derived recursively by (9).

$$\begin{aligned} \boldsymbol{\Psi}(0) &= \mathbf{G}(\boldsymbol{\Phi}(0)), \boldsymbol{\Phi}(1) = \mathbf{F}(\boldsymbol{\Phi}(0), \boldsymbol{\Psi}(0)) \\ \boldsymbol{\Psi}(1) &= \mathbf{G}(\boldsymbol{\Phi}(0), \boldsymbol{\Phi}(1)), \boldsymbol{\Phi}(2) = \mathbf{F}(\boldsymbol{\Phi}(1), \boldsymbol{\Psi}(1)) \\ &\vdots \\ \boldsymbol{\Psi}(k-1) &= \mathbf{G}(\boldsymbol{\Phi}(0) \cdots \boldsymbol{\Phi}(k-1)), \boldsymbol{\Phi}(k) = \mathbf{F}(\boldsymbol{\Phi}(k-1), \boldsymbol{\Psi}(k-1)) \\ \boldsymbol{\Psi}(k) &= \mathbf{G}(\boldsymbol{\Phi}(0) \cdots \boldsymbol{\Phi}(k)) \end{aligned} \quad (9)$$

**Step 3**: Apply the IDT to $\boldsymbol{\Phi}$ and $\boldsymbol{\Psi}$ to obtain the DTM-based SAS in (7). This is equivalent to summating the polynomial terms up to order $K$ with the calculated $\boldsymbol{\Phi}(0)$ to $\boldsymbol{\Phi}(K)$ and $\boldsymbol{\Psi}(0)$ to $\boldsymbol{\Psi}(K)$ as the coefficients based on the definition of IDT in (2).

**Remark 1**: Compared with the process of obtaining an SAS by the ADM, the derivation of an SAS by the DTM is much simpler. The reason is that the DTM determines the symbolic coefficients of SAS terms in a recursive manner for orders from low to high. Thus, those coefficients, i.e. functions of all symbolic variables defined as given in (8) and (9), do not have to explicitly written or stored. When evaluation of the SAS is needed for simulation, those coefficients can easily be calculated from low orders to high orders. In other words, such a recursive procedure to solve and calculate complex symbolic coefficients provides an efficient data structure to generate and store a general SAS for arbitrary orders. That advantage effectively avoids time consuming symbolic computations in derivation of a DTM-based SAS.

**Remark 2**: As for power system models in (5), the coefficients in (7) are computed as follows. Coefficient vectors $\boldsymbol{\Phi}(k) = [\Delta_i(k), W_i(k), E'_{Di}(k), E'_{Qi}(k)]^\mathrm{T} = \mathrm{DT}([\delta_i, \omega_i, e'_{di}, e'_{qi}]^\mathrm{T})$ and $\boldsymbol{\Psi}(k) = [I_{Di}(k), I_{Qi}(k), P_{Ei}(k), E_{Xi}(k), E_{Yi}(k), I_{Xi}(k), I_{Yi}(k), E_{Qi}(k), E_{Qi}(k)]^\mathrm{T} = \mathrm{DT}([i_{di}, i_{qi}, P_{ei}, e_{xi}, e_{yi}, i_{xi}, i_{yi}, e_{di}, e_{qi}]^\mathrm{T})$. Elements of $\boldsymbol{\Phi}(k)$ can be obtained in by applying DT to both sides of (4) with Proposition 1 and Corollary 1. Similarly, Apply DT to both sides of (5) to obtain formulas to calculate elements of $\boldsymbol{\Psi}(k)$.

## III. PROPOSED SEMI-ANALYTICAL SCHEME FOR ONLINE POWER SYSTEM SIMULATION

This section first analyzes the convergence and accuracy for a DTM-based SAS and then designs a multi-time window strategy for transient stability simulation using the SAS. The section also studies how the simulation time cost changes with the order of the SAS for the best time performance. Finally, the section presents the flow chart of a DTM-based semi-analytical scheme for online power system simulation and designs a strategy for paralleling computing.

### A. Convergence and Accuracy of a DTM-Based SAS

The classical SMIB system in (10) is used here to analyze the convergence and accuracy for a DTM-based SAS. The conclusions can be extended to multi-machine systems.

$$\begin{cases} \dot{\delta} = \omega_s \omega \\ \dot{\omega} = \dfrac{1}{2H}\left(P_m - P_{max}\sin\delta - D\omega\right) \end{cases} \quad \begin{cases} \delta(0) = \delta_0 \\ \omega(0) = \omega_0 \end{cases} \quad (10)$$

Here, $H$=3 s, $D$=3 p. u., $\omega_s = 2\pi \times 60$ rad/s, $P_{max}$=1.7 p.u., $P_m$=0.44 p.u. and the initial state has $\delta(0) = 0.26$ rad and $\omega(0) = 0.002$ p.u. Its SAS of order $K$ is given in (11) with coefficients calculated by (12)-(13).

$$\delta(t) = \sum_{k=0}^{K}\Delta(k)\,t^k, \quad \omega(t) = \sum_{k=0}^{K}W(k)\,t^k \quad (11)$$

$$\Delta(k) = \begin{cases} \delta(0), & k = 0 \\ \dfrac{1}{k}\,\omega_s W(k-1), & k \geq 1 \end{cases} \quad (12)$$

$$W(k) = \begin{cases} \omega(0), & k = 0 \\ \dfrac{1}{k}\cdot\dfrac{P_m\boldsymbol{\delta}(k) - P_{max}F(k-1) - DW(k-1)}{2H}, & k \geq 1 \end{cases} \quad (13)$$

The DTM based SAS can be written as either a general recursive form like (11)-(13) or an expanded form for a specific order $K$ similar to an ADM-based SAS. To compare the SAS's from the ADM and DTM and their convergences, Table I lists the coefficients of $t^0$, $t^1$, $t^2$, … in the expressions of the 2nd order and 3rd order SAS's. Here the ADM-based SAS's are from [11] obtained by a modified ADM that gives SAS's as polynomial functions of time for the sake of comparison with polynomial DTM-based SASs. From the table, the differences between two coefficients in row are highlighted. One observation is that all coefficients of a DTM-based SAS remain unchanged when the order increases but the ADM-based SAS continues updating its coefficients. For example, the 2nd order ADM-based SAS lacks "-$D\omega_0$" in the parentheses of the coefficient of $t^2$, which is added when the ADM-based SAS is derived to the 3rd order to become the same as that of the DTM-based SAS. Similarly, to make the coefficient of $t^3$ in the 3rd order SAS be the same with that of the DTM, the ADM needs to derive terms of orders higher than 3. Thus, a DTM-based SAS has better convergence than the ADM-based SAS of the same order because the DTM just needs to derive the $K$-th order SAS to determine the final coefficient of $t^K$ in any SAS's of orders >$K$ while the $K$-th order SAS from the ADM has a less accurate coefficient and has to be extended to a higher order to have the final accurate coefficient.

#### TABLE I
#### COMPARISON BETWEEN THE DTM BASED SAS AND THE ADM BASED SAS

| | | Coefficients of DTM based SAS's | Coefficients of ADM based SAS's |
|---|---|---|---|
| 2nd order SAS | $t^0$ | $\delta_0$ | $\delta_0$ |
| | $t^1$ | $\omega_s\omega_0$ | $\omega_s\omega_0$ |
| | $t^2$ | $\dfrac{\omega_s\omega_b}{4H}(P_m - P_{max}\sin\delta_0 - \boldsymbol{D\omega_0})$ | $\dfrac{\omega_s\omega_b}{4H}(P_m - P_{max}\sin\delta_0)$ |
| 3rd order SAS | $t^0$ | $\delta_0$ | $\delta_0$ |
| | $t^1$ | $\omega_s\omega_0$ | $\omega_s\omega_0$ |
| | $t^2$ | $\dfrac{\omega_s\omega_b}{4H}(P_m - P_{max}\sin\delta_0 - \boldsymbol{D\omega_0})$ | $\dfrac{\omega_s\omega_b}{4H}(P_m - P_{max}\sin\delta_0 - \boldsymbol{D\omega_0})$ |
| | $t^3$ | $-\dfrac{\omega_s^2 P_{max}\omega_0\cos\delta_0}{12H} - \dfrac{\omega_s D\omega_0}{24H^2}(P_m - P_{max}\sin\delta_0 - \boldsymbol{D\omega_0})$ | $-\dfrac{\omega_s^2 P_{max}\omega_0\cos\delta_0}{12H} - \dfrac{\omega_s D\omega_0}{24H^2}(P_m - P_{max}\sin\delta_0)$ |
| … | | … | … |

The rotor angle is calculated by the DTM-based SAS's with order $K$=3, 4, 7 and 15 and then compared with the numerical solution from the RK4 with a time step of 1/1200s as the benchmark. As shown in Fig. 1, the result from each SAS matches the RK4 solution for a certain time window whose length increases with $K$. The maximum time window of accuracy for $K$=3, 4, 7 and 15 are 0.01 s, 0.04 s, 0.10 s and 0.25 s respectively if the error tolerance is set to be $10^{-5}$ rad. Although higher-order SAS's can easily be derived by the DTM for better accuracy, it is observed that for the SMIB system, the time window length saturates when $K$ exceeds 20.

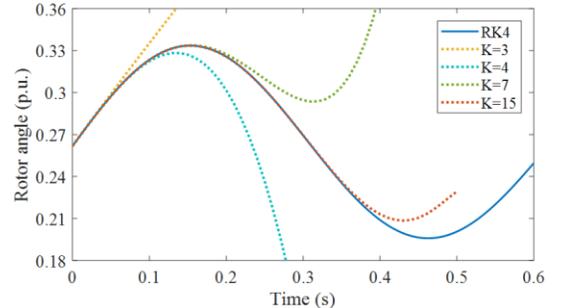

Fig. 1. RK4 solution and the DTM based SAS of different orders.

### B. A Multi-time Window Strategy for the DTM based SAS

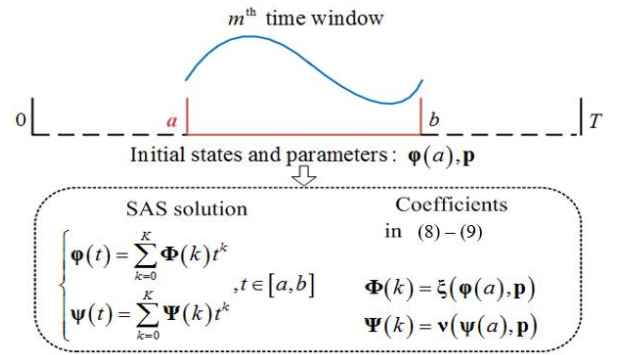

Fig. 2. Multi-time window strategy for proposed method.

Each SAS has a maximum time window to be accurate. Since a simulation period of 5-10 s is usually needed to assess transient stability, a multi-time window strategy is designed to extend the SAS accuracy to the expected simulation period. As shown in Fig. 2, the simulation period is partitioned into a series of time windows, whose length can be pre-determined to meet a desired error tolerance. For the $m$th time window $[a, b]$, the initial state at time $a$ equals the final state of the $(m-1)$th time window. Since coefficients in $\mathbf{\Phi}(k)$ and $\mathbf{\Psi}(k)$ are explicit functions $\mathbf{\xi}$ and $\mathbf{\upsilon}$ of the initial state $\mathbf{\varphi}(a)$ and other symbolized system parameters $\mathbf{p}$, they are directly calculated by (8)-(9) at the beginning of each time window.

## C. Order of an SAS for the Optimal Time Performance

The optimal order of an SAS to be evaluated should minimize the simulation time (i.e. the total time cost for its evaluations over the entire simulation period) with a predefined accuracy. In general, the higher the order of an SAS, the longer time window of accuracy, and vice versa. Note that the simulation time equals the time cost on each evaluation multiplied by the number of evaluations. Let $t_w$ denote the length of the longest time window for expected accuracy. Suppose that the time cost for evaluating the SAS over one time window is $t_{one}$ and the number of time windows to make up simulation period $T$ is $n_w$. Then the total simulation time $t_{total}$ is

$$t_{total} = t_{one} n_w = t_{one} \frac{T}{t_w} \tag{14}$$

When order $K$ rises, time window length $t_w$ is usually increased thus to reduce $n_w$ for the same simulation period $T$. However, the increased order also raises the complexity of the SAS and slows down the evaluation over each window. Therefore, the optimal order should be identified carefully as a tradeoff by checking simulation times using various orders.

## D. Flowchart of the DTM based Online Simulation Scheme

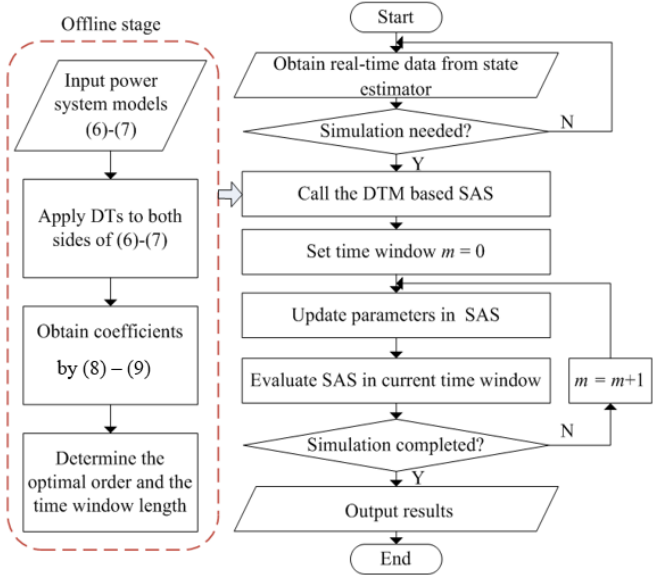

Fig. 3. Flowchart of the DTM-based online simulation scheme.

Fig. 3 shows the flowchart of the DTM based semi-analytical scheme for transient stability simulation. The scheme can be interfaced with the real-time state estimator in the grid control room for online applications. The SAS in (7) is offline derived with coefficients calculated by (8)-(9) in a recursive manner.

Then the prepared SAS together with the optimal order and the time window length is provided for online simulation. Since the SAS explicitly contains the initial state and system parameters as symbolic variables, it can update these values whenever they are changed while its recursive, analytical expression remains the same. For example, when a contingency to be simulated has a fault occur and then cleared, evaluation of the SAS only needs to update values of symbolized elements of the admittance matrix in real time to reflect the topological changes caused by the fault. After evaluations over all time windows are finished, the complete simulation results are obtained by combining trajectories in all time windows.

## E. Evaluating an SAS on Parallel Computers

Applications of HPCs can speed up power system simulation. To demonstrate the superiority of the proposed simulation scheme leveraged by parallel computing, this subsection investigates parallelization of the computations for evaluating a DTM-based SAS.

Like the other types of SAS's [8]-[11], evaluation of a DTM-based SAS naturally fits in a parallel computing environment because it is purely the summation of power series terms in the forms of $\mathbf{\Phi}_i(k)t^k$ and $\mathbf{\Psi}_i(k)t^k$. Further notice that $\mathbf{\Phi}_i(k)t^k$ and $\mathbf{\Psi}_i(k)t^k$ are also summations of several sub-terms shown in (12)-(22). Thus, evaluation of an SAS is the eventually summation of many computation units (CUs), which can be defined as either power series terms or smaller sub-terms depending on the number of available parallel CPU cores.

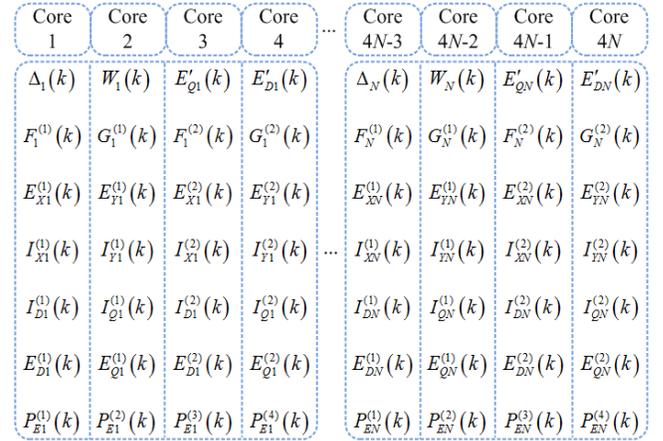

| Core 1 | Core 2 | Core 3 | Core 4 | ... | Core 4N-3 | Core 4N-2 | Core 4N-1 | Core 4N |
|---|---|---|---|---|---|---|---|---|
| $\Delta_1(k)$ | $W_1(k)$ | $E'_{Q1}(k)$ | $E'_{D1}(k)$ | ... | $\Delta_N(k)$ | $W_N(k)$ | $E'_{QN}(k)$ | $E'_{DN}(k)$ |
| $F_1^{(1)}(k)$ | $G_1^{(1)}(k)$ | $F_1^{(2)}(k)$ | $G_1^{(2)}(k)$ | ... | $F_N^{(1)}(k)$ | $G_N^{(1)}(k)$ | $F_N^{(2)}(k)$ | $G_N^{(2)}(k)$ |
| $E_{T1}^{(1)}(k)$ | $E_{T1}^{(1)}(k)$ | $E_{T1}^{(2)}(k)$ | $E_{T1}^{(2)}(k)$ | ... | $E_{TN}^{(1)}(k)$ | $E_{TN}^{(1)}(k)$ | $E_{TN}^{(2)}(k)$ | $E_{TN}^{(2)}(k)$ |
| $I_{X1}^{(1)}(k)$ | $I_{Y1}^{(1)}(k)$ | $I_{X1}^{(2)}(k)$ | $I_{Y1}^{(2)}(k)$ | ... | $I_{XN}^{(1)}(k)$ | $I_{YN}^{(1)}(k)$ | $I_{XN}^{(2)}(k)$ | $I_{YN}^{(2)}(k)$ |
| $I_{D1}^{(1)}(k)$ | $I_{Q1}^{(1)}(k)$ | $I_{D1}^{(2)}(k)$ | $I_{Q1}^{(2)}(k)$ | ... | $I_{DN}^{(1)}(k)$ | $I_{QN}^{(1)}(k)$ | $I_{DN}^{(2)}(k)$ | $I_{QN}^{(2)}(k)$ |
| $E_{D1}^{(1)}(k)$ | $E_{Q1}^{(1)}(k)$ | $E_{D1}^{(2)}(k)$ | $E_{Q1}^{(2)}(k)$ | ... | $E_{DN}^{(1)}(k)$ | $E_{QN}^{(1)}(k)$ | $E_{DN}^{(2)}(k)$ | $E_{QN}^{(2)}(k)$ |
| $P_{E1}^{(1)}(k)$ | $P_{E1}^{(2)}(k)$ | $P_{E1}^{(3)}(k)$ | $P_{E1}^{(4)}(k)$ | ... | $P_{EN}^{(1)}(k)$ | $P_{EN}^{(2)}(k)$ | $P_{EN}^{(3)}(k)$ | $P_{EN}^{(4)}(k)$ |

Fig. 4. Parallel computing in each core.

In this paper, a CU is defined as the DT of a single variable, e.g. $\Delta_i(k)$ and $W_i(k)$. In this context, the number of parallel cores needed is $4N$ for an $N$-machine power system if each generator modeled by a $4^{th}$ order DE. Fig. 4 shows how those $4N$ cores are allocated to evaluate expressions that are either the whole coefficients of the SAS from (8)-(9) or parts of the coefficients. All expressions are reorganized into several groups chronologically. Expressions in each row are independent and calculated by parallel cores while expressions in different rows need to be calculated in a sequential manner. The first row has $4N$ coefficients in $\mathbf{\Phi}(k)$ and needs to be evaluated using all $4N$ cores. Each of the other rows has only 2$N$ or $N$ coefficients and hence partitions each coefficient into 2 or 4 independent expressions for balanced workloads among all cores. The superscripts in each row denotes the independent expressions which can merge to one coefficient by addition.

## IV. CASE STUDY ON IEEE 39-BUS SYSTEM

In this section, a case study on the IEEE 10-machine 39-bus system is conducted to demonstrate the accuracy and time performance of the DTM-based simulation scheme. The simulated contingency is a three-phase fault at the bus 3 cleared after 5 cycles by tripping the line 3 - 4. The RK4 result with a small enough time step of 1/1200 s is used as the benchmark. Both the DTM and RK4 are run in the MATLAB R2017a environment on a laptop with 64-bit Windows operating system and i5-7200U CPU.

### A. Comparing Results from the RK4 and DTM-based SAS

The DTM based SAS trajectory is simulated with the order of $K = 12$ over time windows of 0.2 s, which is 240 times of the RK4 time step. The detailed process for finding the optimal order is presented later in Section IV-B.

Fig. 5 shows the trajectories of rotor angles and speeds from both the DTM-SAS and RK4 including a 1-second pre-fault period. The SAS results accurately match the RK4 results. Fig. 6 further shows the error of the SAS results compared to the the RK4 results as the reference. The maximum error on rotor angles and speeds are within $1.5 \times 10^{-6}$ rad (or p.u.) during the 6-s simulation. As expected, the error reaches the maximum in first few post-fault swings and then significantly decreases.

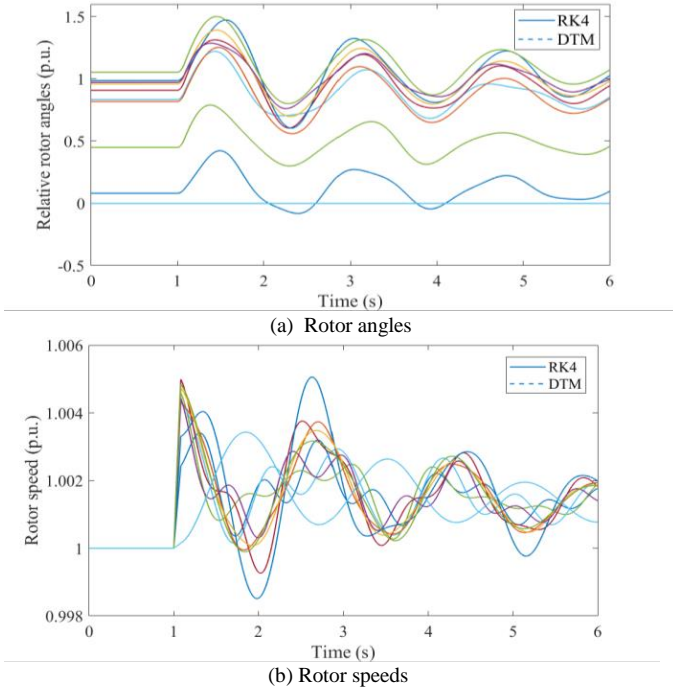

Fig. 5. Comparison of simulation results by RK-4 and the 12th order DTM-based SAS on 0.2 s time windows.

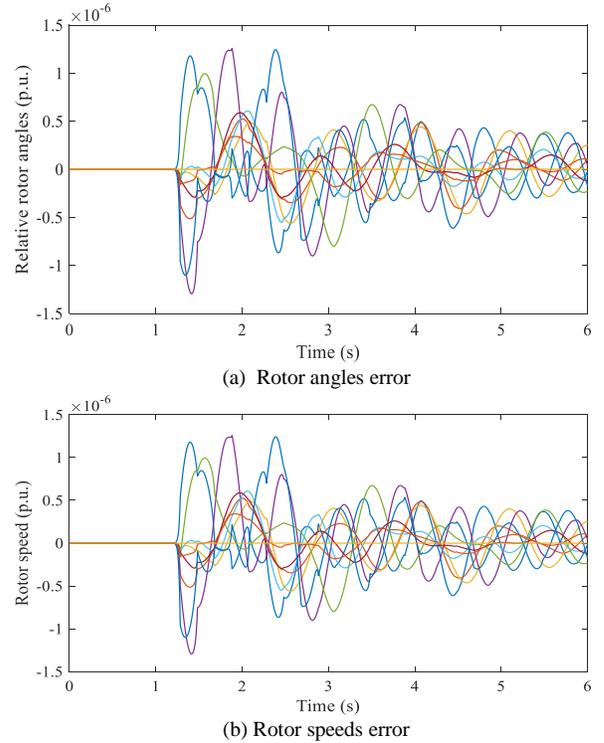

Fig. 6. Error of simulation results by RK-4 and the 12th order DTM-based SAS on 0.2 s time windows.

### B. Order Selection for a DTM based SAS

This section provides a method for determining the optimal order of a DTM-based SAS to minimize simulation time. Table II lists the maximum lengths ($t_w$) of time windows for an SAS with different orders ($K$) to meet different tolerances of accuracy. Obviously, both a bigger tolerance and a higher order enable longer $t_w$. Also, $t_w$ increases rapidly with $K$. When $K$ rises from 6 to 12, $t_w$ is prolonged to 2.4, 4 and 7 times of the original respectively for tolerances of $10^{-3}$, $10^{-5}$ and $10^{-7}$.

TABLE II
THE TIME WINDOW LENGTH UNDER DIFFERENT ORDERS AND ACCURACY

| Tolerance (rad or p.u.) | $K$=6 | $K$=8 | $K$=10 | $K$=12 |
|---|---|---|---|---|
| $10^{-3}$ | 0.12 s | 0.20 s | 0.25 s | 0.30 s |
| $10^{-5}$ | 0.05 s | 0.12 s | 0.16 s | 0.20 s |
| $10^{-7}$ | 0.02 s | 0.06 s | 0.10 s | 0.15 s |

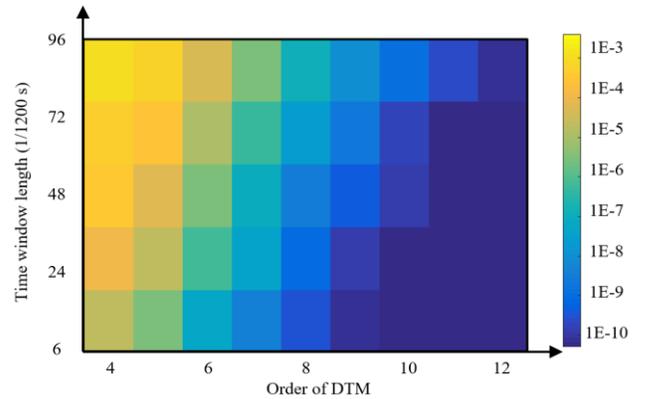

Fig. 7. Relationship among the order, the time window length and error.

Fig. 7 visualizes how the error tolerance is distributed with different combinations on $K$ (horizontal axis) and $t_w$ (vertical axis). $t_w$ in the unit of 1/1200 s indicates how many times it is longer than the time step of the RK4. From the figure, an increase of $K$ or a decrease of $t_w$ has always positive influence on the accuracy. As indicated by the upper-left and bottom-right corners of the figure, the lowest order $K$ paired with the longest time window $t_w$ leads to the largest error while the lowest error occurs at the highest $K$ paired with the smallest $t_w$. Also, the error decreases rapidly with the increase of $K$ when it is not too high. For example, when $K$ increases from 4 to 6 with the time window of 96($\times$1/1200 s), the error is decreased by about one order of magnitude.

Another observation from Fig. 7 is that there are multiple combinations of $K$ and $t_w$ to achieve the same accuracy, as shown in the regions with the same color. Fig. 8 shows by the blue dash line how $t_w$ varies with a change in $K$ to meet tolerance $10^{-5}$. $t_w$ increases rapidly when $K$ grows from 4 to 12 and thereafter, it becomes saturated to reach the longest 0.20 s. The red solid line shows the simulation times by using the SAS's with those orders. As is expected, the minimum simulation time occurs at $K$=12 due to the saturation. From the above analysis, $t_w$ plays a dominant role to speed up simulation when $K$ is low; however, when $K$ becomes high (e.g. 12 for this case), the complexity of the SAS makes its evaluation time over each time window become a dominating factor influencing the overall simulation time. Thus, for this system, the best order is selected as 12 and the corresponding time window $t_w$ is 0.2 s, meaning 30 time windows needed to complete a 6-s simulation by the 12th order SAS.

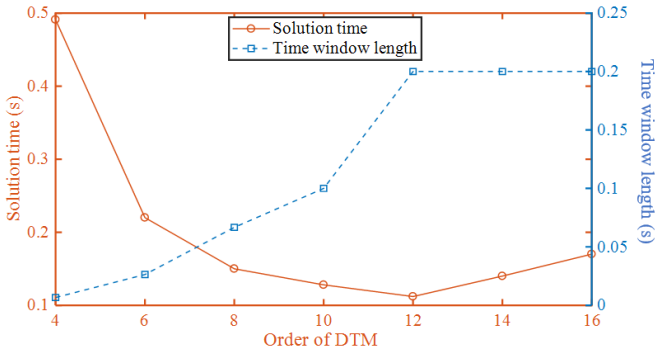

Fig. 8. Relationship among the solution time, the maximum time window length and the DTM order to achieve a pre-defined accuracy of $10^{-5}$ p.u..

### C. Time Performance Comparison

This section demonstrates the time performance of a DTM based SAS compared with numerical integration by the RK4. The CPU time for the RK4 is 1.207 s in MATLAB as the benchmark for comparison. A speed up index is calculated by dividing the RK4 time cost by the SAS-based simulation time equaling to the total SAS evaluation time over the desired simulation period with or without parallel computing.

Tab III lists the values of the SAS-based simulation time and the speed up index to achieve different accuracies using the 12th order SAS. The values before and after "/" are respectively for scenarios without and with parallelization. From the table, the SAS gains more than 10 times speed up for tolerance $10^{-5}$ without parallelization of any computation, and the speed up index drops to 7.5 when tolerance $10^{-9}$ is required. When ideal parallelization on at least 40 cores is considered that ignores time costs on communications among cores, the parallelization strategy proposed in Section III-E expedites the SAS-based simulation to be 300-400 times faster. The main contributor to the speed up index with the SAS-based simulation is the largely reduced number of time windows compared to the time step of the RK4 for the same accuracy. Fig. 9 indicates the time points (i.e. starting points of time windows) when the SAS is evaluated to calculate the rotor angle of generator 8. Only 30 time windows are needed to complete the 6-s simulation and each SAS time window has 240 RK4 time steps.

TABLE III
Time Performances on the SAS with/without Parallelization

| Tolerance | Time cost on SAS-based simulation (s) | Speed up Index |
|---|---|---|
| $10^{-3}$ | 0.110 / 0.00275 * | 11.0 / 439 * |
| $10^{-5}$ | 0.118 / 0.00295 | 10.2 / 409 |
| $10^{-7}$ | 0.146 / 0.00365 | 8.3 / 331 |
| $10^{-9}$ | 0.160 / 0.00401 | 7.5 / 301 |

* Values before and after "/" are without and with parallelizations, respectively.

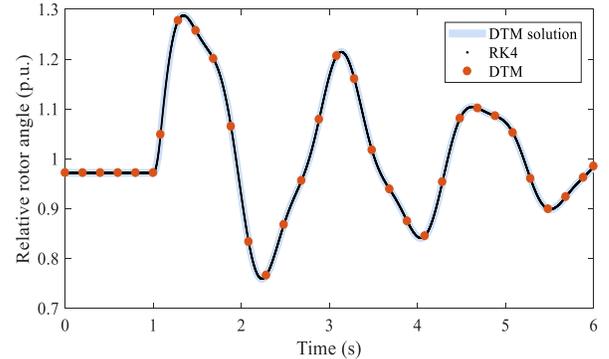

Fig. 9. Comparison of the number of time windows given by RK-4 and the 12th order DTM-based SAS.

## V. Conclusion

This paper presents a DTM-based semi-analytical approach for transient stability simulation of multi-machine power systems. Compared with an SAS derived from the ADM, a DTM-based SAS adopts a general recursive form for any order $K$ to avoid complex symbolic computations and can also be evaluated to a higher order. A multi-time window strategy is then designed to guarantee the DTM-based SAS accuracy over a series of time windows to make up the desired simulation period. Since the same level of accuracy can be achieved by multiple combinations of orders and time window lengths regarding a DTM-based, the optimal order is studied to minimize simulation time. A DTM-based semi-analytical scheme for online simulation is also suggested. Parallelization of the SAS-based simulation has been studied to demonstrate the superiority of DTM-based SAS's leveraged by parallel computing. Test results on IEEE 39-bus system have shown that a DTM-based SAS can considerably speed up transient stability simulation while keeping high accuracy because it allows a high order SAS compared to an ADM-based SAS.

## REFERENCES

[1] S. Zadkhast, J. Jatskevich, E. Vaahedi, "A multi-decomposition approach for accelerated time-domain simulation of transient stability problems," *IEEE Trans. on Power Syst.*, vol. 30, no. 5, pp. 2301-2311, Sep. 2015.

[2] Y. Liu, Q. Jiang, "Two-stage parallel waveform relaxation method for large-scale power system transient stability simulation," *IEEE Trans. on Power Syst.*, vol. 31, no. 1, pp. 153-162, Jan. 2016.

[3] P. Aristidou, S. Lebeau, T. V. Cutsem, "Power system dynamic simulations using a parallel two-level Schur-complement decomposition," *IEEE Trans. on Power Syst.*, vol. 31, no. 5, pp. 3984-3995, Sep. 2016.

[4] G. Gurrala, A. Dimitrovski, P. Sreekanth, S. Simunovic, M. Starke, "Parareal in time for fast power system dynamic simulations," *IEEE Trans. Power Syst.*, vol. 31, no. 3, pp. 1820-1830, Jul. 2016.

[5] N. Duan, A. Dimitrovski, S. Simunovic, K. Sun, "Applying reduced generator models in the coarse solver of parareal in time parallel power system simulation," in Proc. *IEEE PES Innovative Smart Grid Technologies Conference Europe*, 2016, pp. 1-5.

[6] R. Diao, S. Jin, F. Howell, "On parallelizing single dynamic simulation using HPC techniques and APIs of commercial software," *IEEE Trans. on Power Syst.*, vol. 32, no. 3, pp. 2225-2233, May. 2017.

[7] I. Konstantelos, G. Jamgotchian, S. H. Tindemans, P. Duchesne, S. Cole, C. Merckx, G. Strbac, P. Panciatici, "Implementation of a massively parallel dynamic security assessment platform for large-scale grids," *IEEE Trans. on Smart Grid*, vol. 8, no. 3, pp. 1417-1426, May 2017.

[8] E. Abreut, B. Wang, and K. Sun, "Semi-analytical fault-on trajectory simulation and its application in direct methods," in Proc. *Power and Energy Society General Meeting*, 2017, pp. 1-5.

[9] N. Duan and K. Sun, "Finding semi-analytic solutions of power system differential-algebraic equations for fast transient stability simulation," Dec. 2014. [Online]. Available: https://arxiv.org/abs/1412.0904.

[10] N. Duan, K. Sun, "Power system simulation using the multi-stage Adomian Decomposition Method," *IEEE Trans. on Power Syst.*, vol. 32, no. 1, pp 430-441, Jan. 2017.

[11] N. Duan, K. Sun, "Application of the Adomian Decomposition Method for semi-analytic solutions of power system differential algebraic equations," in Proc. *IEEE PowerTech*, 2015, pp. 1-5.

[12] C. Liu, B. Wang, K. Sun, "Fast power system simulation using semi analytical solutions based on Pade Approximants," in Proc. *Power and Energy Society General Meeting*, 2017, pp. 1-5.

[13] G. Gurrala, D. L. Dinesha, A. Dimitrovski, P. Sreekanth, S. Simunovic, M. Starke, "Large multi-machine power system simulations using multi-Stage Adomian Decomposition," *IEEE Trans. on Power Syst.*, vol. 32, No. 5, pp. 3594-3606, January 2017.

[14] G. Gurrala, A. Dimitrovski, P. Sreekanth, S. Simunovic, M. Starke, K. Sun, "Application of Adomian Decomposition for Multi-Machine Power System Simulation," 2015 IEEE PES General Meeting.

[15] B. Wang, N. Duan, K. Sun, "A Time-Power Series Based Semi-Analytical Approach for Power System Simulation," *IEEE Transactions on Power Systems*, submitted.

[16] E. Pukhov, G. Georgii, "Differential transforms and circuit theory," *International Journal of Circuit Theory and Applications*, vol. 10, no. 3, pp. 265-276, Jun. 1982.

[17] I. H. A. Hassan, "Application to differential transformation method for solving systems of differential equations," *Applied Mathematical Modelling*, vol. 32, no. 12, pp. 2552-2559, Oct. 2007.

[18] M. Ghafarian, A. Ariaei, "Free vibration analysis of a system of elastically interconnected rotating tapered Timoshenko beams using differential transform method," *International Journal of Mechanical Sciences*, vol. 107, pp. 93-109, 2016.

[19] M. Matinfar, S. R. Bahar, M. Ghasemi, "Solving the Lienard equation by differential transform method," *World Journal of Modelling and Simulation*, vol. 8, no. 2, pp. 142-146, 2012.

[20] L. Xie, C. Zhou, S. Xu, "An effective numerical method to solve a class of nonlinear singular boundary value problems using improved differential transform method," *SpringerPlus*, vol. 5, no. 1, pp. 1066, Jul. 2016.

[21] F. Kenmogne, "Generalizing of differential transform method for solving nonlinear differential equations," *Applied & Computational Mathematics*, vol. 4, no. 1, pp. 1000196, Jan. 2015.

[22] V. S. Erturk, Z. M. Odibat, S. Momani, "The multi-step differential transform method and its application to determine the solutions of non-linear oscillators," *Advances in Applied Mathematics and Mechanics*, vol. 4, no. 4, pp. 422-438, Aug. 2012.

[23] C. W. Berta, H. Zeng, "Analysis of axial vibration of compound bars by differential transformation method," *Journal of Sound and Vibration*, vol. 275, no. 3, pp. 641-647, Aug. 2004.

[24] Z. Odibat, S. Momani, "A generalized differential transform method for linear partial differential equations of fractional order," Applied Mathematics Letters, vol. 21, no. 2, pp. 194-199, Feb. 2008.